# An example of a braided category of $C^*$-algebras[*]

Małgorzata Rowicka - Kudlicka
Dept. of Math. Methods in Physics, Univ. of Warsaw,
Hoża 69, 00-682 Warszawa, Poland
Dept. of Mathematics, Univ. of Munich,
Theresienstr. 39, 80333 München, Germany
rowicka@fuw.edu.pl

February 12, 1997, revised November 5, 1997

**Abstract**

A braided category of $C^*$-algebras is constructed. Its objects are $C^*$-algebras endowed with an action of the group $\mathbb{R}$, its morphisms are $C^*$-algebra morphisms intertwining the action of $\mathbb{R}$, the braided tensor product of its two objects essentially depends on the action of $\mathbb{R}$ on considered $C^*$-algebras. Braided tensor products of any object with $C_\infty(\mathbb{R})$ and $C(\mathbb{T}^1)$ and $C(\mathbb{T}^2)$ are discussed.

## 0 Introduction

As far as the usual *tensor product* $\otimes$ of algebras A and B is concerned, the subalgebras $A \otimes I$ and $I \otimes B$ commute, i.e. for any $a \in A$ and $b \in B$,

$$j_1(a)j_2(b) = j_2(b)j_1(a)$$

where $j_1$ is an injection $A \hookrightarrow A \otimes B$, such that $j_1(a) = a \otimes I$ and $j_2$ is an injection $B \hookrightarrow A \otimes B$, such that $j_2(b) = I \otimes b$.

The simplest nontrivial commutation rule occurs while a *supertensor product* of superalgebras (i.e. $\mathbb{Z}_2$-graded algebras) A and B is considered:

$$j_1(a)j_2(b) = (-1)^{|a||b|} j_2(b)j_1(a)$$

where $|a|, |b|$ are the grades of homogenous elements $a \in A$ and $b \in B$ respectively.

A further generalization provides a notion of a *braided tensor product* $\underline{\otimes}$ of algebras. Majid in his book [4, Chapter 9] showed that the natural setting for such products is a braided tensor category. He proposed the following commutation rule:

$$j_1(a)j_2(b) = \sum_k j_2(b_k)j_1(a_k)$$

---

[*]Partly supported by Deutscher Akademischer Austauschdienst

where $a_k \in A$, $b_k \in B$ and the sum is finite. The commutation rule can be encoded in a mapping $\Psi_{A,B} : j_1(A)j_2(B) \to j_2(B)j_1(A)$, such that $\Psi_{A,B}(j_1(a)j_2(b)) = \sum j_2(b_k)j_1(a_k)$. Mapping $\Psi_{A,B}$ is called braiding since it can be thought of as braids (generating the Artin braid group). Precisely, although both $\Psi_{A,B}$ and $(\Psi_{B,A})^{-1}$ are mappings from $j_1(A)j_2(B)$ into $j_2(B)j_1(A)$ they are in general distinct (see below).

$$\Psi_{A,B} = \begin{array}{c} j_1(a) \quad j_2(b) \\ \diagdown\!\!\!\!\diagup \\ \end{array} \neq \begin{array}{c} j_1(a) \quad j_2(b) \\ \diagup\!\!\!\!\diagdown \\ \end{array} = (\Psi_{B,A})^{-1}$$

Woronowicz in [7] cast Majid's braided tensor product of algebras into a $C^*$-algebra setting. Woronowicz introduced a "crossed product", or what we will call throughout the paper a *braided tensor product* $\otimes_\mathcal{C}$ *of $C^*$-algebras* ($C^*$-braided tensor product for short) as a main concept. He claimed that Majid's approach based on braiding mappings given in advance is not general enough, since the braiding $\Psi$ may have a too small domain to contain all the information required. Namely, in Woronowicz's axioms both of the sets $j_1(A)j_2(B)$ and $j_2(B)j_1(A)$ are dense in a $C^*$-algebra $A \otimes_\mathcal{C} B$, but they can have a trivial intersection.

Within this setting, Woronowicz built in [7] a braided category of $C^*$-algebras endowed with an action of the group $\mathbb{Z} \times \mathbb{S}^1$.

In this paper, following [7], we formally build another braided category of $C^*$-algebras, precisely a *braided category of $C^*$-algebras endowed with an action of the group $\mathbb{R}$*.

Although the group $\mathbb{R}$ appears easier to handle than the group $\mathbb{Z} \times \mathbb{S}^1$, the construction of a braided category of $C^*$-algebras with an action of $\mathbb{R}$ is more difficult.

In section 1 we recall Woronowicz's construction and proceed with constructing a braided tensor product algebra of $C^*$-algebras with an action of $\mathbb{R}$. Presented proofs are different from those in [7]. Most difficult is proving that a $C^*$-braided tensor product is a $C^*$-algebra. To this end one has to discover and investigate the commutation rule involved. It turns out that it is enough to know that for elements $a_n$, $b_n$ from a certain dense subset of $A \otimes_\mathcal{C} B$:

$$j_2(b_n)j_1(a_n) = C_1(n) \int_{\mathbb{R}^2} j_1(\alpha_r^A(a))j_2(\alpha_t^B(b))e^{-C_2(n)(r^2+t^2)-C_3(n)irt} dr dt$$

where

- $a_n, b_n$ are constructed from $a \in A$ and $b \in B$ respectively by "smashing" it.
- $\alpha^A, \alpha^B$ are actions of $\mathbb{R}$ on $C^*$-algebras $A$ and $B$, respectively
- $C_1(n), C_2(n), C_3(n)$ are certain functions of $n \in \mathbb{Z}$

(see Lemmas 1.4 - 1.5 for details).

Then we construct a $C^*$-braided category with an action of $\mathbb{R}$. We used freely the propositions proved by Woronowicz in [7], they are true also in our case. The difference between groups $\mathbb{R}$ and $\mathbb{S}^1 \times \mathbb{Z}$ matters only when assiociativity of the $\otimes_\mathcal{C}$ functor is proved. Theorem 2.7 in section 2 is the main result of this paper. We prove therein that one can freely make $C^*$– braided tensor product among a collection of $C^*$-algebras endowed with an action of $\mathbb{R}$ in a consistent way.

In section 3 we discuss the $C^*$-braided tensor products of any object of our category with $C^*$-algebras $C_\infty(\mathbb{R})$ and $C(\mathbb{T}^1)$ and $C(\mathbb{T}^2)$. We show that they are in fact isomorphic with more familiar crossed products of a $C^*$-algebra by locally compact abelian groups $\mathbb{R}$, $\mathbb{Z}$ and $\mathbb{Z} \times \mathbb{Z}$, respectively.

## 0.1 Notation

We follow the notation of [7]. Let $H$ be a separable Hilbert space and let $C^*(H)$ denote the set of all separable nondegenerate $C^*$-subalgebras of $B(H)$ (a $C^*$-algebra $A$ is nondegerate when $AH$ is dense in $H$).

Let $A \in C^*(H)$. The multiplier algebra of a $A$ is denoted by $\mathrm{M}(A)$ and

$$\mathrm{M}(A) = \{a \in B(H): \ ab, \ ba \in A \text{ for any } b \in A\}$$

$\mathrm{M}(A)$ is a unital $C^*$-subalgebra of $B(H)$ and $A$ is an ideal in $\mathrm{M}(A)$. We endow $\mathrm{M}(A)$ with the strict topology, i.e. the weakest topology such that for all $a \in A$ the mappings

$$\mathrm{M}(A) \ni x \mapsto ax \in A \quad \mathrm{M}(A) \ni x \mapsto xa \in A$$

are continuous. The strict topology is weaker than the norm topology.

Let $A$ be a separable $C^*$-algebra and $B \in C^*(H)$. The set of all morphisms acting from $A$ into $B$, will be denoted by $\mathrm{Mor}(A, B)$. We recall that $\varphi \in \mathrm{Mor}(A, B)$ if $\varphi$ is a *-algebra morphism mapping $A$ into $\mathrm{M}(B)$, such that $\varphi(A)B$ is dense in $B$. Any morphism $\varphi \in \mathrm{Mor}(A, B)$ admits a unique extension to a unital *-algebra homomorphism from $M(A)$ into $M(B)$. We endow $\mathrm{Mor}(A, B)$ with the weakest topology such that for all mappings

$$\mathrm{Mor}(A, B) \ni \varphi \mapsto \varphi(a) \in \mathrm{M}(B)$$

are strictly continuous.

We use the minimal tensor product of $C^*$-algebras. Precisely, for any $A \in C^*(H)$ and $B \in C^*(K)$ the tensor product $A \otimes B$ is the closure of $A \otimes_{alg} B$ in the operator norm from $B(H \otimes K)$ (see e.g. [5, 8.15.15]). Of course, $A \otimes B$ is a $C^*$-algebra.

Let $I_A$ denote the unit of $\mathrm{M}(A)$ and $\mathrm{id}_A$ be the identity morphism on $A$.

For any locally compact topological space $\mathbb{X}$, $C_\infty(\mathbb{X})$ will denote the $C^*$-algebra of all continuous functions vanishing at infinity on $\mathbb{X}$. Then $\mathrm{M}(C_\infty(\mathbb{X})) = C_\mathrm{b}(\mathbb{X})$ is the algebra of all continuous bounded functions on $\mathbb{X}$. The algebra of all compact operators acting on Hilbert space $H$ will be denoted by $CB(H)$. $\mathrm{M}(CB(H))$ coincides with the algebra $B(H)$ of all bounded operators acting on $H$. From now on we write $K$ instead of $CB(L^2(\mathbb{R}))$.

$C^*$-algebras endowed with an action of the group $\mathbb{R}$ will play a crucial role later on. For our purpose we need a more abstract definition of an action than the usual one. Let $e \in \mathrm{Mor}(C_\infty(\mathbb{R}), \mathbb{C})$ and $\Delta_\mathbb{R} \in \mathrm{Mor}(C_\infty(\mathbb{R}), C_\infty(\mathbb{R}) \otimes C_\infty(\mathbb{R}))$ be the standard counit and comultiplication, i.e. $e(f) = f(0)$ and $(\Delta_\mathbb{R} f)(x, y) = f(x + y)$ for any $f \in C_\infty(\mathbb{R})$ and $x, y \in \mathbb{R}$. *An action of the group $\mathbb{R}$ on an $C^*$-algebra $A$ is a morphism* $\alpha^A \in \mathrm{Mor}(A, A \otimes C_\infty(\mathbb{R}))$ *such that* $(\mathrm{id} \otimes e)\alpha^A = \mathrm{id}$ and the diagram

$$\begin{CD}
A @>{\alpha^A}>> A \otimes C_\infty(\mathbb{R}) \\
@V{\alpha^A}VV @VV{\alpha^A \otimes \mathrm{id}}V \\
A \otimes C_\infty(\mathbb{R}) @>>{\mathrm{id} \otimes \Delta_\mathbb{R}}> A \otimes C_\infty(\mathbb{R}) \otimes C_\infty(\mathbb{R})
\end{CD}$$

commutes.

Let $A$ be a $C^*$-algebra and $\alpha^A$ be an action of the group $\mathbb{R}$ on it. Hence we obtain a family of automorphisms $(\alpha_t^A)_{t \in \mathbb{R}}$ of $A$:

$$\alpha_t^A = (\mathrm{id} \otimes \chi_t)\alpha^A \tag{1}$$

where the evaluation $\chi_t \in \mathrm{Mor}(C_\infty(\mathbb{R}), \mathbb{C})$ is given by $\chi_t(f) = f(t)$. Moreover

$$\alpha_0 = \mathrm{id}_A \quad \text{and} \quad \forall t, s \in \mathbb{R} \quad \alpha_t \alpha_s = \alpha_{t+s} ,$$

the latter follows from the commutativity of the preceding diagram. Let $T_r$ denote the translation operator defined for $f \in C_\mathrm{b}(\mathbb{R})$ and $r \in \mathbb{R}$ by $(T_r f)(x) := f(x + r)$. Then $T_r = (\mathrm{id} \otimes \chi_r)\Delta_\mathbb{R}$

# 1 Braided tensor product of $C^*$-algebras

We introduce the definition of $C^*$-braided tensor product from [7].[1] It is a generalization of the usual $C^*$-tensor product. It means, that if we set $j_1 \in \mathrm{Mor}(A, A \otimes B)$, such that $j_1(a) = a \otimes I_B$ and $j_2 \in \mathrm{Mor}(B, A \otimes B)$ such that $j_2(b) = I_A \otimes b$ then $(j_1, j_2, A \otimes B)$ will satisfy the below definition.

**Definition 1.1** *Let $A,B,C$ be $C^*$-algebras, $j_1 \in \mathrm{Mor}(A, C)$, $j_2 \in \mathrm{Mor}(B, C)$ and $j : A \otimes_{alg} B \to \mathrm{M}(C)$ be a linear map, such that:*

$$j(\sum_{k=1}^n a_k \otimes b_k) = \sum_{k=1}^n j_1(a_k) j_2(b_k) \tag{2}$$

*for any $a_1, a_2, ..., a_n \in A, b_1, b_2, ..., b_n \in B$.*

*We say that $(j_1, j_2, C)$ is a braided tensor product of $C^*$-algebras $A$ and $B$ if*

1. *$j(A \otimes_{alg} B)$ is a dense subset of $C$,*
2. *For any $\xi \in A \otimes_{alg} B$, $\{j(\xi) = 0\} \Rightarrow \{\xi = 0\}$.*

**Remark 1.2** *Let $A$, $B$ and $D$ be $C^*$-algebras and $(j_1, j_2, C)$ be a $C^*$-braided tensor product of $A$ and $B$. Then any morphism $\varphi \in \mathrm{Mor}(C, D)$ is uniquely determined by $\varphi \circ j_1$ and $\varphi \circ j_2$ ([7]).*

---

[1] there, the name *a crossed product of $C^*$-algebras* is used instead

We construct an example of the $C^*$-braided tensor product similar to that presented in [7], but with the group $\mathbb{R}$ in place of $\mathbb{S}^1 \times \mathbb{Z}$.

Let $\hat{q}$ and $\hat{p}$ be position and momentum operators acting on $L^2(\mathbb{R})$, i.e. $\hat{q}$ denotes the multiplication by the coordinate $q$ and $\hat{p} = \frac{1}{i}\frac{d}{dq}$. We shall use the functional calculus of self–adjoint operators. We define linear mappings $i_1, i_2 : C_\infty(\mathbb{R}) \to B(L^2(\mathbb{R}))$ by

$$i_1(f) = f(\hat{q}) \quad \text{and} \quad i_2(g) = g(\hat{p})$$

for any $f, g \in C_\infty(\mathbb{R})$. We see that $i_1, i_2 \in \text{Mor}(C_\infty(\mathbb{R}), K)$, so there is a unique extension of $i_1, i_2$ to maps $C_b(\mathbb{R}) \to B(L^2(\mathbb{R}))$.

Here and subsequently in this section, $A, B$ are $C^*$-algebras, endowed with the action of $\mathbb{R}$ on it, which is denoted by $\alpha^A, \alpha^B$ respectively. We define an automorphism $\sigma_{B,A} \in \text{Mor}(B \otimes A, A \otimes B)$ by

$$\sigma_{B,A}(b \otimes a) = a \otimes b$$

For any $a \in A$ and $b \in B$ we set

$$\begin{aligned} j_1(a) &= (\sigma_{B,A} \otimes i_1)(I_B \otimes \alpha^A(a)), \\ j_2(b) &= I_A \otimes (\text{id}_B \otimes i_2)\alpha^B(b). \end{aligned} \quad (3)$$

It is easy to check that $j_1 \in \text{Mor}(A, A \otimes B \otimes K)$ and $j_2 \in \text{Mor}(B, A \otimes B \otimes K)$.

Let $j : A \otimes_{alg} B \to M(A \otimes B \otimes K)$ be the linear map introduced by (2) and

$$A \otimes_\mathcal{C} B = \overline{j(A \otimes_{alg} B)}^{\|\cdot\|} \quad (4)$$

**Theorem 1.3** *A triple $(j_1, j_2, A \otimes_\mathcal{C} B)$ is a $C^*$-braided tensor product of $A$ and $B$.*

A scheme of our proof:

$$\left.\begin{array}{l} \text{Lem. 1.4} \implies \text{Lem. 1.5} \implies \text{Lem. 1.6} \implies \text{Prop. 1.7} \\ \phantom{\text{Lem. 1.4} \implies \text{Lem. 1.5} \implies} \text{Lem. 1.8} \implies \text{Prop. 1.9} \end{array}\right\} \implies \text{Th. 1.3}$$

The main difficulty is to prove Prop. 1.7, i.e. that $A \otimes_\mathcal{C} B$ is a $C^*$-algebra. To this end we use Lemmas 1.4 – 1.6. The basic trick is to work with a dense subset of elements, for which the commutation rule within $A \otimes_\mathcal{C} B$ is relatively easy. Lemma 1.4 states that the difficulty lies in interchanging elements in the third leg $(K)$ of $A \otimes_\mathcal{C} B \subset A \otimes B \otimes K$. In Lemma 1.5 we prove that Lemma 1.4 is in fact all we need to interchange elements from the dense subset of $A \otimes_\mathcal{C} B$ mentioned above. Lemma 1.6 is just a last step toward Prop.1.7. In proving Proposition 1.9 the way we act is similar. To begin with, we solve pur problem only in the third leg. Next we prove that that in fact all the work is already done.

Let $G_n$ (n=1,2,...) be a $\delta$ - like sequence consisting of normalized Gauss functions $G_n(r) := \sqrt{\frac{n}{2\pi}} e^{-\frac{n}{2}r^2}$ and let $\mathcal{G}_n(r,t) := \frac{n}{2\pi\sqrt{n^2+1}} \exp\left\{-\frac{\frac{n}{2}(r^2+t^2)+n^2 irt}{n^2+1}\right\}$. One can verify that

$$\int_{\mathbb{R}^2} \mathcal{G}_n(r-q, t-p) e^{ip(q-\tilde{q})} dp = G_n(r-\tilde{q}) e^{-\frac{(q-\tilde{q})^2}{2n}} e^{it(q-\tilde{q})} \quad (5)$$

**Lemma 1.4** Let $f, g \in C_{\mathrm{b}}(\mathbb{R})$,

$$f_n(x) = \int_{\mathbb{R}} (T_r f)(x) \mathrm{G}_n(r) dr \quad \text{and} \quad g_n(x) = \int_{\mathbb{R}} (T_r g)(x) \mathrm{G}_n(r) dr. \tag{6}$$

Then

$$i_2(g_n) i_1(f_n) = \int_{\mathbb{R}^2} i_1(T_r f) i_2(T_t g) \mathcal{G}_n(r, t) dr dt. \tag{7}$$

**Proof:** Let $\mathcal{F}$ denote the Fourier transform: for any $f \in L^1(\mathbb{R})$ $(\mathcal{F}f)(p) := \frac{1}{\sqrt{2\pi}} \int_{\mathbb{R}} f(q) e^{-ipq} dq$. There is a unique extension of $\mathcal{F}$ to a unitary operator acting on $L^2(\mathbb{R})$. It is known that $\hat{p} = \mathcal{F}^{-1} \hat{q} \mathcal{F}$. Therefore for any $f \in C_{\mathrm{b}}(\mathbb{R})$ $i_2(f) = f(\hat{p}) = \mathcal{F}^{-1} f(\hat{q}) \mathcal{F} = \mathcal{F}^{-1} i_1(f) \mathcal{F}$. Let $\mathcal{S}(\mathbb{R})$ be the Schwartz space of test functions and $\psi, \varphi \in \mathcal{S}(\mathbb{R})$. Let us compute matrix elements of left-hand side of (7):

$$(\psi | LHS | \varphi) = (\psi | i_2(g_n) i_1(f_n) | \varphi) = (\psi | \mathcal{F}^{-1} \mathcal{F} i_2(g_n) \mathcal{F}^{-1} \mathcal{F} i_1(f_n) | \varphi)$$

Since for any $f \in C_{\mathrm{b}}(\mathbb{R})$, $\psi \in \mathcal{S}(\mathbb{R})$ we have $i_2(f)\psi = \mathcal{F}^{-1} i_1(f) \mathcal{F} \psi$. It implies that

$$(\psi | LHS | \varphi) = (\mathcal{F}\psi | i_1(g_n) \mathcal{F} i_1(f_n) | \varphi)$$

$$= \frac{1}{\sqrt{2\pi}} \int_{\mathbb{R}} \overline{(\mathcal{F}\psi)}(p)(g_n)(p) \int_{\mathbb{R}} f_n(q) e^{-ipq} \varphi(q) dp dq =$$

$$= \int_{\mathbb{R}^2} g(t) f(r) \left( \frac{1}{\sqrt{2\pi}} \int_{\mathbb{R}^2} \overline{(\mathcal{F}\psi)(p)} e^{-ipq} \varphi(q) \mathrm{G}_n(t-p) \mathrm{G}_n(r-q) dp dq \right) dr dt =$$

$$= \int_{\mathbb{R}^2} f(r) g(t) d\mu_{\psi, \varphi, n}(r, t)$$

where

$$d\mu_{\psi, \varphi, n}(r, t) = \left( \frac{1}{\sqrt{2\pi}} \int_{\mathbb{R}^2} \overline{(\mathcal{F}\psi)(p)} e^{-ipq} \varphi(q) \mathrm{G}_n(t-p) \mathrm{G}_n(r-q) dp dq \right) dr dt$$

Let us observe that the measure density

$$\frac{d\mu_{\psi, \varphi, n}(r, t)}{dr dt} = \frac{1}{\sqrt{2\pi}} \int_{\mathbb{R}^2} \overline{\left( \int_{\mathbb{R}} \psi(q') e^{-ipq'} dq' \right)} e^{-ipq} \varphi(q) \mathrm{G}_n(t-p) \mathrm{G}_n(r-q) dp dq =$$

$$= \frac{1}{\sqrt{2\pi}} \int_{\mathbb{R}^2} \overline{\psi(q')} \varphi(q) \mathrm{G}_n(r-q) \left( \int_{\mathbb{R}} \mathrm{G}_n(t-p) e^{ip(q'-q)} dp \right) dq' dq =$$

$$= \frac{1}{\sqrt{2\pi}} \int_{\mathbb{R}^2} \overline{\psi(q')} \varphi(q) \mathrm{G}_n(r-q) e^{-\frac{(q'-q)^2}{2n}} e^{it(q'-q)} dq' dq$$

Similarly, for the right-hand side of (7) we have

$$(\psi | RHS | \varphi) = (\psi | \int_{\mathbb{R}^2} i_1(T_r f) i_2(T_t g) \mathcal{G}_n(r, t) dr dt | \varphi) =$$

$$= \int_{\mathbb{R}^2} (\psi | i_1(T_r f)\mathcal{F}^{-1}\mathcal{F} i_2(T_t g)\mathcal{F}^{-1}|\mathcal{F}\varphi) \mathcal{G}_n(r,t) dr dt =$$

$$= \frac{1}{\sqrt{2\pi}} \int_{\mathbb{R}^2} dr dt \mathcal{G}_n(r,t) \int_{\mathbb{R}^2} \overline{\psi(q)} f(q+r) g(p+t) e^{ipq} (\mathcal{F}\varphi)(p) dq dp =$$

$$= \int_{\mathbb{R}^2} f(r) g(t) \left( \frac{1}{\sqrt{2\pi}} \int_{\mathbb{R}^2} \mathcal{G}_n(r-q, t-p) \overline{\psi(q)} e^{ipq} (\mathcal{F}\varphi)(p) dq dp \right) dr dt =$$

$$= \int_{\mathbb{R}^2} f(r) g(t) d\nu_{\psi,\varphi,n}(r,t)$$

where

$$d\nu_{\psi,\varphi,n} = \left( \frac{1}{\sqrt{2\pi}} \int_{\mathbb{R}^2} \mathcal{G}_n(r-q, t-p) \overline{\psi(q)} e^{ipq} (\mathcal{F}\varphi)(p) dq dp \right) dr dt$$

By (5), the measure density

$$\frac{d\nu_{\psi,\varphi,n}(r,t)}{dr dt} = \frac{1}{\sqrt{2\pi}} \int_{\mathbb{R}^2} \mathcal{G}_n(r-q, t-p) \overline{\psi(q)} e^{ipq} \left( \int_{\mathbb{R}} \varphi(\tilde{q}) e^{-ip\tilde{q}} d\tilde{q} \right) dq dp =$$

$$= \frac{1}{\sqrt{2\pi}} \int_{\mathbb{R}^2} \overline{\psi(q)} \varphi(\tilde{q}) \left( \int_{\mathbb{R}} \mathcal{G}_n(r-q, t-p) e^{ip(q-\tilde{q})} dp \right) dq d\tilde{q} =$$

$$= \frac{1}{\sqrt{2\pi}} \int_{\mathbb{R}^2} \overline{\psi(q)} \varphi(\tilde{q}) G_n(r-\tilde{q}) e^{-\frac{(q-\tilde{q})^2}{2n}} e^{it(q-\tilde{q})} dq d\tilde{q} = \frac{d\mu_{\psi,\varphi,n}(r,t)}{dr dt}$$

And (7) follows. $\square$

**Lemma 1.5** *Let $a \in A$ and $b \in B$,*

$$a_n = \int \alpha_r^A(a) G_n(r) dr \quad \text{and} \quad b_n = \int \alpha_t^B(b) G_n(t) dt. \tag{8}$$

*Then*

$$j_2(b_n) j_1(a_n) = \int_{\mathbb{R}^2} j_1(\alpha_r^A(a)) j_2(\alpha_t^B(b)) \mathcal{G}_n(r,t) dr dt. \tag{9}$$

**Proof:** Let $\omega$ and $\rho$ be continuous linear functionals on $C^*$-algebras $A$ and $B$ respectively. Define $f, g, f_n, g_n \in C_{\mathrm{b}}(\mathbb{R})$ by $f(t) = (\omega \otimes \mathrm{id})(\alpha_t^A(a))$ and $g(t) = (\rho \otimes \mathrm{id})(\alpha_t^B(b))$ and $f_n(t) = (\omega \otimes \mathrm{id})(\alpha_t^A(a_n))$ and $g_n(t) = (\rho \otimes \mathrm{id})(\alpha_t^B(b_n))$. It is easily seen that

$$(\omega \otimes \rho \otimes \mathrm{id}) j_2(b_n) j_1(a_n) = (\rho \otimes i_2)\alpha^B(b_n)(\omega \otimes i_1)\alpha^A(a_n)$$

We will show that $f$ and $f_n$ as well as $g$ and $g_n$ are related by (6). Observe that

$$(\omega \otimes i_1)\alpha^A(a_n) = i_1 \left( \int (\omega \otimes \mathrm{id}) \alpha^A (\mathrm{id} \otimes \chi_r) \alpha^A(a) G_n(r) dr \right)$$

Using (1) and following remarks we obtain

$$(\omega \otimes \mathrm{id}) \alpha^A (\mathrm{id} \otimes \chi_r) \alpha^A(a) = (\omega \otimes \mathrm{id} \otimes \chi_r)(\alpha^A \otimes \mathrm{id}) \alpha^A(a) =$$

$$= (\omega \otimes \mathrm{id} \otimes \chi_r)(\mathrm{id} \otimes \Delta_{\mathbb{R}})\alpha^A(a) =$$
$$= T_r\{(\omega \otimes \mathrm{id})\alpha^A(a)\} = T_r(f)$$

Therefore
$$(\omega \otimes i_1)\alpha^A(a_n) = i_1(\int_{\mathbb{R}} T_r(f)\mathrm{G}_n(r)dr)$$

In the similar manner we get
$$(\rho \otimes i_2)\alpha^B(b_n) = i_2(\int_{\mathbb{R}} T_r(g)\mathrm{G}_n(r)dr)$$

On the other hand we conclude:
$$(\omega \otimes \rho \otimes \mathrm{id})j_1(\alpha_r^A(a))j_2(\alpha_t^B(b)) = (\omega \otimes i_1)\alpha^A(\alpha_r^A(a))(\rho \otimes i_2)\alpha^B(\alpha_t^B(b)) = i_1(T_r f)i_2(T_t g)$$

Hence
$$(\omega \otimes \rho \otimes \mathrm{id})\int_{\mathbb{R}^2} j_1(\alpha_r^A(a)) = \int_{\mathbb{R}^2} i_1(T_r f)i_2(T_t g)\mathcal{G}_n(r,t)drdt$$

Now Lemma 1.4 makes (9) obvious. □

**Lemma 1.6** *For any $a \in A$ and $b \in B$, $j_2(b)j_1(a) \in A \otimes_{\mathcal{C}} B$.*

**Proof:** According to Definition (1.1), $j_1\left(\alpha_r^A(a)\right) j_2\left(\alpha_t^B(b)\right) = \mathrm{j}\left(\alpha_r^A(a) \otimes \alpha_t^B(b)\right) \in A \otimes_{\mathcal{C}} B$ for any $r, t \in \mathbb{R}$. Let $a_n \in A$ and $b_n \in B$ be the elements introduced by (8). By (9) we see that $j_2(b_n)j_1(a_n) \in A \otimes_{\mathcal{C}} B$. Remembering that $\mathrm{G}_n$ is a sequence of $\delta$-like functions and $\alpha_r^A(a)$ and $\alpha_t^B(b)$ are norm continuous with respect to $r, t \in \mathbb{R}$ ([7]) we conclude that $a_n \to a$ and $b_n \to b$ in norm topology. It shows that $j_2(b)j_1(a) = \lim_{n\to\infty} j_2(b_n)j_1(a_n) \in A \otimes_{\mathcal{C}} B$. □

**Proposition 1.7** $A \otimes_{\mathcal{C}} B$ *is a $C^*$-algebra, $j_1 \in \mathrm{Mor}(A, A \otimes_{\mathcal{C}} B)$, $j_2 \in \mathrm{Mor}(B, A \otimes_{\mathcal{C}} B)$.*

**Proof:** According to definition of $A \otimes_{\mathcal{C}} B$
$$j_1(A)(A \otimes_{\mathcal{C}} B) \subset A \otimes_{\mathcal{C}} B \quad \text{and} \quad (A \otimes_{\mathcal{C}} B)j_2(B) \subset A \otimes_{\mathcal{C}} B$$

Applying Lemma 1.6 we get $j_1(c)j_2(b)j_1(a)j_2(d) \in A \otimes_{\mathcal{C}} B$ for any $a, c \in A$ and $b, d \in B$. Hence $A \otimes_{\mathcal{C}} B$ is closed under multiplication. Again by (1.6) $(j_1(a)j_2(b))^* = j_2(b^*)j_1(a^*) \in A \otimes_{\mathcal{C}} B$. Thus we proved that $A \otimes_{\mathcal{C}} B$ is a $C^*$-algebra.

By Lemma 1.6 we see that
$$(A \otimes_{\mathcal{C}} B)j_1(A) \subset A \otimes_{\mathcal{C}} B \quad \text{and} \quad j_2(B)A \otimes_{\mathcal{C}} B \subset A \otimes_{\mathcal{C}} B$$

Hence $j_1(A), j_2(B) \subset \mathrm{M}(A \otimes_{\mathcal{C}} B)$. Since $j_1(A)$ is a $C^*$-algebra $j_1(A)j_1(A)$ is dense in $j_1(A)$, hence $j_1(A)j_1(A)j_2(B)$ is dense in $j_1(A)j_2(B)$, and $j_1(A)j_2(B)$ is dense in $A \otimes_{\mathcal{C}} B$ by definition. The analogous conclusion can be drawn for $j_2$. We have thus proved that $j_1 \in \mathrm{Mor}(A, A \otimes_{\mathcal{C}} B)$ and $j_2 \in \mathrm{Mor}(B, A \otimes_{\mathcal{C}} B)$. □

**Lemma 1.8** *Let $f_1, f_2, ..., f_n, g_1, g_2, ..., g_n \in C_b(\mathbb{R})$.*
*Assume that $\sum_{k=1}^n i_1(f_k)i_2(g_k) = 0$. Then $\sum_{k=1}^n f_k \otimes g_k = 0$.*

**Proof:** We may assume that $f_k$ are linearly independent. For any $\psi \in L^2(\mathbb{R})$ we have

$$(\sum_{k=1}^n f_k(\hat{q})g_k(\hat{p})\psi)(x) = 0.$$

Inserting in this formula $T_r\psi$ instead of $\psi$ and remembering that the translation operators commute with $\hat{p}$ we get

$$\left(\sum_{k=1}^n i_1(f_k)i_2(g_k)T_r\psi\right)(x) = \sum_{k=1}^n f_k(x)\left(g_k(\hat{p})\psi\right)(x+r) = 0$$

for almost all x and all r. Introducing new variables we see that

$$\sum_{k=1}^n f_k(x)(g_k(\hat{p})\psi)(y) = 0$$

for almost all $(x, y) \in \mathbb{R}^2$. Since functions $f_k$ are linearly independent it follows that $i_2(g_k)\psi = 0$ for $k = 1, 2, ..., n$ and for any $\psi \in L^2(\mathbb{R})$. Hence $g_k = 0$ for $k = 1, 2, ..., n$ and therefore

$$\sum_{k=1}^n f_k \otimes g_k = 0$$

□

**Proposition 1.9** $\forall \xi \in A \otimes_{alg} B, \{j(\xi) = 0\} \Rightarrow \{\xi = 0\}$.

**Proof:** Let $\xi = \sum_{k=1}^n a_k \otimes b_k$ where $a_k \in A$ and $b_k \in B$. Assume that $j(\xi) = 0$. Then for any continuous linear functionals $\omega$ on $A$ and $\rho$ on $B$:

$$0 = (\omega \otimes \rho \otimes \text{id})j(\sum_{k=1}^n a_k \otimes b_k) = \sum_{k=1}^n (\omega \otimes i_1)\alpha^A(a_k)(\rho \otimes i_2)\alpha^B(b_k) = \sum_{k=1}^n i_1(f_k) \otimes i_2(g_k)$$

where $f_k, g_k \in C_b(\mathbb{R})$ are given by $f_k = (\omega \otimes \text{id})\alpha^A(a_k)$ and $f_k = (\rho \otimes \text{id})\alpha^B(b_k)$. By Lemma 1.8 we get $\sum_{k=1}^n f_k \otimes g_k = 0$. Hence $\left(\sum_{k=1}^n f_k \otimes g_k\right)(0, 0) = 0$. On the other hand, $\left(\sum_{k=1}^n f_k \otimes g_k\right)(t, \tau) = \sum_{k=1}^n \omega(\alpha_t^A(a_k))\rho(\alpha_\tau^B(b_k))$. Therefore:

$$0 = \sum_{k=1}^n \omega(\alpha_0^A(a_k))\rho(\alpha_0^B(b_k)) = (\omega \otimes \rho)\sum_{k=1}^n a_k \otimes b_k = (\omega \otimes \rho)(\xi)$$

Remembering that $\omega$ and $\rho$ were arbitrarily chosen we conclude that $\xi = 0$. □

Armed with Prop. 1.7 and Prop. 1.9 we are ready for:

**Proof of Theorem 1.3:** By Proposition 1.7 $A \otimes_C B$ is a $C^*$-subalgebra of $M(A \otimes B \otimes K)$, $j_1 \in \text{Mor}(A, A \otimes_C B)$, $j_2 \in \text{Mor}(B, A \otimes_C B)$ and by Proposition 1.9 the second condition of Definition 1.1 is fullfilled. Since the first condition is fullfilled automatically in our case (see (4)), we proved that $(j_1, j_2, A \otimes_C B)$ is a $C^*$-braided tensor product of $A$ and $B$. □

# 2 Braided categories of $C^*$-algebras

In the first subsection we introduce necessary definitions. In the second one we built, out of $C^*$-braided tensor product considered in section 1, an example of a braided category of $C^*$-algebras endowed with actions of $\mathbb{R}$.

## 2.1 Introduction

We assume that the reader is familiar with the basic notions of category theory, such as category itself, its objects and its morphisms. We also use notions of functors and natural projections between two categories. Details can be found e. g. in [4, Subsection 9.1], which is also a good starting point to problems considered in this paper.

We present Woronowicz's definitions of $C^*$– monoidal category and $C^*$-braided category, which are $C^*$-algebras' counterparts of well-known notions of monoidal category (see [3]) and a braided monoidal category ( see [4, Subsection 9.2]), respectively.

Here we present a definition of a $C^*$-monoidal category from [7]. The class of objects of our category $\mathcal{C}$ is denoted by $Ob_\mathcal{C}$. For any $A, B \in Ob_\mathcal{C}$ let $\mathrm{Mor}_\mathcal{C}(A, B)$ denote a set of morphisms acting from $A$ into $B$. Let $\mathrm{Proj}_1$ and $\mathrm{Proj}_2$ denote the canonical projections acting from $\mathcal{C} \times \mathcal{C}$ onto $\mathcal{C}$.

**Definition 2.1** *Let $\mathcal{C}$ be a category and let $\otimes_\mathcal{C}$ be a covariant functor acting from $\mathcal{C} \times \mathcal{C}$ into $\mathcal{C}$. Let $j_1$ and $j_2$ be natural maps acting from $\mathrm{Proj}_1$ and $\mathrm{Proj}_2$ respectively into $\otimes_\mathcal{C}$. We say that $(\mathcal{C}, \otimes_\mathcal{C}, j_1, j_2)$ is a $C^*$-**monoidal category** if*
*1. for any $A, B, C \in Ob_\mathcal{C}$ and any $\varphi, \varphi' \in \mathrm{Mor}_\mathcal{C}(A \otimes_\mathcal{C} B, C)$*

$$\left\{ \begin{array}{l} \varphi \circ j_1(A, B) = \varphi' \circ j_1(A, B) \\ \varphi \circ j_2(A, B) = \varphi' \circ j_2(A, B) \end{array} \right\} \Rightarrow \{\varphi = \varphi'\}$$

*2. $\otimes_\mathcal{C}$ is associative, i.e. for any $A, B, C \in Ob_\mathcal{C}$, there exists an isomorphism*

$$\Psi_{ABC} \in \mathrm{Mor}_\mathcal{C}((A \otimes_\mathcal{C} B) \otimes_\mathcal{C} C, A \otimes_\mathcal{C} (B \otimes_\mathcal{C} C))$$

*such that*

$$\left. \begin{array}{rcl} j_1(A, B \otimes_\mathcal{C} C) &=& \Psi_{ABC} \circ j_1(A \otimes_\mathcal{C} B, C) \circ j_1(A, B) \\ j_2(A, B \otimes_\mathcal{C} C) \circ j_1(B, C) &=& \Psi_{ABC} \circ j_1(A \otimes_\mathcal{C} B, C) \circ j_2(A, B) \\ j_2(A, B \otimes_\mathcal{C} C) \circ j_2(B, C) &=& \Psi_{ABC} \circ j_2(A \otimes_\mathcal{C} B, C) \end{array} \right\} \qquad (10)$$

In the above definition $\otimes_\mathcal{C}$ is a binary operation acting on objects and morphisms of $\mathcal{C}$. Mapping $j_1$ is a natural mapping acting from $\mathrm{Proj}_1$ into $\otimes_\mathcal{C}$, i.e. for any $A, B \in Ob_\mathcal{C}$, $j_1(A, B) \in \mathrm{Mor}_\mathcal{C}(A, B)$. Moreover, for any $\varphi \in \mathrm{Mor}_\mathcal{C}(A, A')$ and $\psi \in \mathrm{Mor}_\mathcal{C}(B, B')$ the diagram

$$
\begin{array}{ccc}
A & \xrightarrow{j_1(A,B)} & A \otimes_{\mathcal{C}} B \\
\varphi \downarrow & & \downarrow \varphi \otimes_{\mathcal{C}} \psi \\
A' & \xrightarrow{j_1(A',B')} & A' \otimes_{\mathcal{C}} B'
\end{array}
$$

commutes (analogously for $j_2(A,B)$). Condition I means that $A \otimes_{\mathcal{C}} B$ is in a certain sense generated by $j_1(A)$ and $j_2(B)$. Therefore isomorphism $\Psi_{ABC}$ satisfying condition II is unique. For more details, we refer the reader to [7, Section 2]. We follow [7] in defining a braided category of $C^*$-algebras.

**Definition 2.2** *Let $(\mathcal{C}, \otimes_{\mathcal{C}}, j_1, j_2)$ be a $C^*$-monoidal category such that objects of $\mathcal{C}$ are $C^*$-algebras endowed with an additional structure and, for any $A, B \in Ob_{\mathcal{C}}$, $\mathrm{Mor}_{\mathcal{C}}(A, B)$ is a subset of $\mathrm{Mor}(A, B)$ consisting of all morphisms preserving this additional structure. We call $(\mathcal{C}, \otimes_{\mathcal{C}}, j_1, j_2)$ a braided category of $C^*$-algebras if for any $A, B \in Ob_{\mathcal{C}}$, $(j_1(A,B), j_2(A,B), A \otimes_{\mathcal{C}} B)$ is a $C^*$-braided tensor product of $A$ and $B$.*

## 2.2 An example of $C^*$-braided category endowed with actions of $\mathbb{R}$

The previously introduced braided tensor product of $C^*$-algebras endowed with actions of group $\mathbb{R}$ gives rise to an interesting $C^*$-braided category $(\mathcal{C}, \otimes_{\mathcal{C}}, j_1, j_2)$. Objects of $\mathcal{C}$ are $C^*$-algebras endowed with actions of $\mathbb{R}$, $\otimes_{\mathcal{C}}$ is the $C^*$-braided tensor product introduced in Definition 1.1. For any $A, B \in Ob_{\mathcal{C}}$, $j_1(A,B) \in \mathrm{Mor}(A, A \otimes_{\mathcal{C}} B)$, $j_2(A,B) \in \mathrm{Mor}(B, A \otimes_{\mathcal{C}} B)$ coincide with the morphisms $j_1, j_2$ introduced by (3) in Section 1. Morphisms of the category $\mathcal{C}$ are $C^*$– morphisms intertwining the actions of $\mathbb{R}$: for any $A, B \in Ob_{\mathcal{C}}$,

$$\mathrm{Mor}_{\mathcal{C}}(A, B) = \{\varphi \in \mathrm{Mor}(A, B) : \alpha^B \circ \varphi = (\varphi \otimes \mathrm{id}_{C_{infty}(\mathbb{R})}) \circ \alpha^A\} \tag{11}$$

At the moment, the $C^*$-algebras $A \otimes_{\mathcal{C}} B$ $(A, B \in Ob_{\mathcal{C}})$ are not endowed with an action of $\mathbb{R}$ and consequently $j_1(A,B)$ and $j_2(A,B)$ are not $\mathcal{C}$-morphisms yet. To achieve the construction of $(\mathcal{C}, \otimes_{\mathcal{C}}, j_1, j_2)$ we have to introduce a $\otimes_{\mathcal{C}}$-product of $\mathcal{C}$- morphisms and a natural action of $\mathbb{R}$ on all $A \otimes_{\mathcal{C}} B$ and then to show that $j_1(A,B), j_2(A,B)$ and the $\otimes_{\mathcal{C}}$-product of $\mathcal{C}$-morphisms are $\mathcal{C}$-morphisms.

Repeating arguments used in [7] ( Propositions 2.3-2.4 and further remarks) we may prove the Propositions 2.3 – 2.4 and Remarks 2.5 – 2.6.

**Proposition 2.3** *Let any $A, A', B, B' \in Ob_{\mathcal{C}}$, $\varphi \in \mathrm{Mor}_{\mathcal{C}}(A, A')$ and $\psi \in \mathrm{Mor}_{\mathcal{C}}(B, B')$. Then there exists a unique $\varphi \otimes_{\mathcal{C}} \psi \in \mathrm{Mor}\,(A \otimes_{\mathcal{C}} B, \mathrm{M}(A' \otimes_{\mathcal{C}} B'))$ such that*

$$(\varphi \otimes_{\mathcal{C}} \psi) \circ j_1(A,B) = j_1(A',B') \circ \varphi,$$
$$(\varphi \otimes_{\mathcal{C}} \psi) \circ j_2(A,B) = j_2(A',B') \circ \psi,$$
$$\varphi \otimes_{\mathcal{C}} \psi \in \mathrm{Mor}(A \otimes_{\mathcal{C}} B, A' \otimes_{\mathcal{C}} B')$$

If moreover $A'', B'' \in \text{Ob}_\mathcal{C}$, $\varphi' \in \text{Mor}_\mathcal{C}(A', A'')$ and $\psi' \in \text{Mor}_\mathcal{C}(B', B'')$, then

$$(\varphi' \otimes_\mathcal{C} \psi') \circ (\varphi \otimes_\mathcal{C} \psi) = \varphi'\varphi \otimes_\mathcal{C} \psi'\psi$$

**Proposition 2.4** *For any $A, B \in \text{Ob}_\mathcal{C}$ there exists a unique action $\alpha^{A \otimes_\mathcal{C} B}$ of $\mathbb{R}$ on the $C^*$-braided tensor product $A \otimes_\mathcal{C} B$ such that*

$$\alpha_t^{A \otimes_\mathcal{C} B} = \alpha_t^A \otimes_\mathcal{C} \alpha_t^B$$

*for any $t \in \mathbb{R}$. Thus, for any $A, B \in \text{Ob}_\mathcal{C}$, $A \otimes_\mathcal{C} B$ is an object of the category $\mathcal{C}$.*

**Remark 2.5** *Let the assumptions of Proposition 2.3 hold. Then*

$$\varphi \otimes_\mathcal{C} \psi \in \text{Mor}_\mathcal{C}(A \otimes_\mathcal{C} B, A' \otimes_\mathcal{C} B')$$

**Remark 2.6** *Mappings $j_1$ and $j_2$, given by (3), are natural maps from respectively $\text{Proj}_1$ and $\text{Proj}_2$ into $\otimes_\mathcal{C}$.*

The main result of this paper is

**Theorem 2.7** $(\mathcal{C}, \otimes_\mathcal{C}, j_1, j_2)$ *is a braided category of $C^*$-algebras.*

**Proof:** We will follow [7]. First we prove that $(\mathcal{C}, \otimes_\mathcal{C}, j_1, j_2)$ is a $C^*$-monoidal category. Observe that for any $A, B \in \text{Ob}_\mathcal{C}$ morphisms $j_1(A, B)$ and $j_2(A, B)$ are natural maps from functors $\text{Proj}_1$ and $\text{Proj}_2$ respectively into $\otimes_\mathcal{C}$. Condition 1 of definition 2.1 is automatically fullfilled by Remark 1.2. What is left is to show that $\otimes_\mathcal{C}$ is associative.

Let us introduce an operator $W$ acting on $L^2(\mathbb{R}) \otimes L^2(\mathbb{R})$

$$\forall f \in L^2(\mathbb{R}) \quad (Wf)(x, y) = f(y - x, x)$$

Note that $W$ is a unitary operator.

It is easy to check that

$$\left.\begin{array}{rcl} W(\hat{q} \otimes I + I \otimes \hat{q}) & = & (I \otimes \hat{q})W \\ W(\hat{p} \otimes I + I \otimes \hat{q}) & = & (I \otimes \hat{p} + \hat{q} \otimes I)W \\ W(I \otimes \hat{p}) & = & (\hat{p} \otimes I + I \otimes \hat{p})W \end{array}\right\} \quad (12)$$

For any bounded operator $Q$ acting on $L^2(\mathbb{R}) \otimes L^2(\mathbb{R})$ we set

$$\psi(Q) = WQW^*$$

Then $\psi$ is an isomorphism acting on $K \otimes K$.

We claim that every $f \in C_\infty(\mathbb{R})$ satisfies

$$\left.\begin{array}{rcl} \psi\{(i_1 \otimes i_1) \circ \Delta_\mathbb{R}(f)\} & = & I_K \otimes i_1(f) \\ \psi\{(i_2 \otimes i_1) \circ \Delta_\mathbb{R}(f)\} & = & (i_1 \otimes i_2) \circ \Delta_\mathbb{R}(f) \\ \psi(I_K \otimes i_2(f)) & = & (i_2 \otimes i_2) \circ \Delta_\mathbb{R}(f) \end{array}\right\} \quad (13)$$

From (12) it follows that (13) holds for the function $x$, where $x(t) = t$ for any $t \in \mathbb{R}$. Since $C_\infty(\mathbb{R})$ is generated by $x$ ([8]), we conclude that (13) is true for any $f \in C_\infty(\mathbb{R})$.

Let
$$\Psi_{ABC} : A \otimes B \otimes K \otimes C \otimes K \to A \otimes B \otimes C \otimes K \otimes K$$
be a map given by
$$\Psi_{ABC} = \mathrm{id}_{A \otimes B} \otimes [(\mathrm{id}_C \otimes \psi) \circ (\sigma_{K,C} \otimes \mathrm{id}_K)] \tag{14}$$

Since $\psi|_{K \otimes K}$ is an isomorphism, then so is $\Psi_{ABC}$. Using properties of an action of $\mathbb{R}$, (3), (11) and (14) we see that (10) coincides with (13). It remains to prove that
$$\Psi_{ABC} \in \mathrm{Mor}_\mathcal{C} \left( (A \otimes_\mathcal{C} B) \otimes_\mathcal{C} C, A \otimes_\mathcal{C} (B \otimes_\mathcal{C} C) \right).$$

Remembering that for any $C^*$-algebras $A$ and $B$, $j_1(A)j_2(B)$ is dense in $A \otimes_\mathcal{C} B$, we see that $\Psi_{ABC}$ maps a dense subset of $(A \otimes_\mathcal{C} B) \otimes_\mathcal{C} C$ onto a dense subset of $A \otimes_\mathcal{C} (B \otimes_\mathcal{C} C)$.

Therefore
$$\Psi_{ABC} \in \mathrm{Mor}_\mathcal{C}((A \otimes_\mathcal{C} B) \otimes_\mathcal{C} C, A \otimes_\mathcal{C} (B \otimes_\mathcal{C} C))$$
□

## 3 Examples

To illustrate the above definitions let us consider the following three distinguished objects of the described category $\mathcal{C}$: $C_\infty(\mathbb{R})$ and $C(\mathbb{T}^1)$ and $C(\mathbb{T}^2)$, (to shorten notation we write $\mathbb{T}^1$ instead of $\mathbb{R}/\frac{2\pi}{\mathcal{T}}\mathbb{Z}$ and $\mathbb{T}^2$ instead of $\mathbb{R}/\frac{2\pi}{\mathcal{T}_1}\mathbb{Z} \times \mathbb{R}/\frac{2\pi}{\mathcal{T}_2}\mathbb{Z}$ where $\mathcal{T}, \mathcal{T}_1, \mathcal{T}_2 \in \mathbb{R}_+$ and $\mathcal{T}_1/\mathcal{T}_2$ is not rational).

Observe that there is a map $\rho$:
$$\rho : \mathbb{R} \ni t \mapsto (t + \frac{2\pi}{\mathcal{T}_1}\mathbb{Z}, t + \frac{2\pi}{\mathcal{T}_2}\mathbb{Z}) \in \mathbb{T}^2$$

Therefore $C(\mathbb{T}^2)$ can be treated also as a subalgebra of $C_\mathrm{b}(\mathbb{R})$.

Obviously $C_\infty(\mathbb{R}) \subset C_\mathrm{b}(\mathbb{R})$ and $C(\mathbb{T}^1) \subset C_\mathrm{b}(\mathbb{R})$.

The group $\mathbb{R}$ acts on these three $C^*$-algebras by translations. Therefore for any object $A \in \mathcal{C}$ we can form $C^*$-braided tensor products $A \otimes_\mathcal{C} C_\infty(\mathbb{R})$ and $A \otimes_\mathcal{C} C(\mathbb{T}^1)$ and $A \otimes_\mathcal{C} C(\mathbb{T}^2)$. We will show that these $C^*$-braided tensor products are isomorphic to the crossed products of $A$ with groups $\mathbb{R}$ and $\mathcal{T}\mathbb{Z}$ and $\mathcal{T}_1\mathbb{Z} \times \mathcal{T}_2\mathbb{Z}$ respectively in the sense explained below. Crossed products of a $C^*$-algebra with a locally compact group acting on it were investigated, among others, by Landstad ([2]). Let us recall their definition ([1] Section 2.7.1.,[2]).

Let $A$ be a $C^*$-algebra endowed with an action $\alpha$ of a locally compact abelian group $G$. In $C_\mathrm{o}(G, A)$ (the norm continuous $A$ - valued functions on $G$ with compact support) we define a convolution product, an involution and a norm by

$$(\varphi\theta)(g) = \int_G \varphi(h)\alpha_{-h}\theta(g-h)dh,$$

$$\varphi^*(g) = \alpha_{-g}(\varphi(-g)^*),$$

$$||\varphi||_1 = \int_G ||\varphi(g)||dg,$$

for any $\varphi, \theta \in C_o(G, A)$. Thus $C_o(G, A)$ becomes a normed $*$-algebra. Now, we introduce a new norm on $C_o(G, A)$

$$||\varphi|| = \sup_\pi ||\pi(\varphi)||,$$

where $\pi$ ranges over all the Hilbert space representations of $C_o(G, A)$. One can show ([1] Section 2.7.1.) that $||\cdot||$ is a $C^*-$ norm. The completion of $C_o(G, A)$ in this norm is called a $C^*$- *crossed product of A and G* and is denoted by $A \otimes_\alpha G$.

If $A \in C^*(H)$, $A \in \mathcal{C}$ and $G \subset \mathbb{R}$, then we may consider $A \otimes_\alpha G$ as an element of $C^*(L^2(\mathbb{R}, H))$; the action of $\varphi \in C_o(G, A) \subset A \otimes_\alpha G$ on a vector $v \in L^2(\mathbb{R}, H)$ is given by

$$(\varphi v)(t) = \int_G \alpha_t(\varphi(g))v(t-g)dg \quad t \in \mathbb{R}, \ g \in G$$

We can assume that $\varphi(g) = k(g)a$, where $k \in C_o(G)$ and $a \in A$. Then

$$(\varphi v)(t) = \int_G \alpha_t(a)k(g)v(t-g)dg = \alpha_t(a)\int_G k(g)v(t-g)dg \tag{15}$$

It is known ([2]) that $C^*-$ crossed product of $A$ with an abelian $G$ is spanned by operators $\varphi$ where $\varphi \in C_o(G, A)$.

First we shall consider $A \otimes_\mathcal{C} C_\infty(\mathbb{R})$. In this case $G = \mathbb{R}$ For any $k \in L^1(\mathbb{R})$ we define the Fourier transform by:

$$(\mathcal{F}k)(p) := \frac{1}{\sqrt{2\pi}}\int_\mathbb{R} k(x)e^{-ixp}dx$$

Notice that $\mathcal{F}k$ is a continuous function on $\mathbb{R}$ tending to zero at infinity: $\mathcal{F}k \in C_\infty(\mathbb{R})$. Then (15) takes form

$$(\varphi v)(t) = \alpha_t(a)\int_\mathbb{R} k(s)v(t-s)ds =$$

$$= \alpha_t(a)\int_\mathbb{R} k(s)(e^{-is\hat{p}}v)(t)ds =$$

$$= \alpha_t(a)\{(\mathcal{F}k)(\hat{p})v\}(t) =$$

$$= [(\mathrm{id}_A \otimes i_1)\alpha(a)\{I_A \otimes (i_2(\mathcal{F}k)v)\}](t)$$

Therefore

$$\varphi = (\mathrm{id}_A \otimes i_1)\alpha(a)\{I_A \otimes i_2(\mathcal{F}k)\}$$

where $\mathcal{F}k$ runs over a dense subset of $C_\infty(\mathbb{R})$ ([6] Theorem 1.2.4.).

Hence $(\mathrm{id}_A \otimes i_1)\alpha^A(A)(I_A \otimes i_2(C_\infty(\mathbb{R})))$ is dense in $A \otimes_\alpha \mathbb{R}$.

On the other hand ([8], Formula 2.5)

$$A \otimes_\mathcal{C} C_\infty(\mathbb{R}) \subset M\left(A \otimes C_\infty(\mathbb{R}) \otimes K\right) = C_b(\mathbb{R}, \mathrm{M}(A \otimes K))$$

where $C_{\mathrm{b}}(\mathbb{R}, \mathrm{M}(A \otimes K))$ is a set of all bounded strictly continuous $\mathrm{M}(A \otimes K)$ - valued functions on $\mathbb{R}$. For any $\xi \in A \otimes_{\mathcal{C}} C_{\infty}(\mathbb{R})$ and $t \in \mathbb{R}$, let $\xi[t] \in \mathrm{M}(A \otimes K)$ denote the value of the function $\xi \in C_{\mathrm{b}}(\mathbb{R}, \mathrm{M}(A \otimes K))$ at the point $t \in \mathbb{R}$:

$$\xi[t] = (\mathrm{id}_A \otimes \chi_t \otimes \mathrm{id}_K)\xi$$

where $\chi_t \in \mathrm{Mor}(C_{\infty}(\mathbb{R}), \mathbb{C})$ is the evaluation functional. With this notation definition (3) can be rewritten as follows

$$(j_1(a))[t] = (\mathrm{id} \otimes i_1)\alpha^A(a)$$

and

$$(j_2(b))[t] = (I \otimes (\mathrm{id} \otimes i_2)\alpha^B(b))[t] = I \otimes b(t+\hat{p}) = I \otimes e^{-it\hat{q}} b(\hat{p}) e^{it\hat{q}} =$$
$$= (I \otimes e^{-it\hat{q}})(I \otimes i_2(b))(I \otimes e^{it\hat{q}})$$

Hence, observing that $j_1(a)$ and $e^{it\hat{q}}$ commute, we obtain

$$(j_1(a)j_2(b))[t] = (I \otimes e^{-it\hat{q}})(\mathrm{id} \otimes i_1)\alpha^A(a)(I \otimes i_2(b))(I \otimes e^{it\hat{q}})$$

and therefore $(\mathrm{id} \otimes i_1)\alpha^A(A)(I \otimes i_2(C_{\infty}(\mathbb{R})))$ is unitarily equivalent to $j_1(A)j_2(C_{\infty}(\mathbb{R}))$, which is (see Definition 1.1) dense in $A \otimes_{\mathcal{C}} C_{\infty}(\mathbb{R})$. Thus we see that $C^*$-algebras $A \otimes_{\mathcal{C}} C_{\infty}(\mathbb{R})$ and $A \otimes_{\alpha} \mathbb{R}$ are isomorphic

$$A \otimes_{\mathcal{C}} C_{\infty}(\mathbb{R}) \simeq A \otimes_{\alpha} \mathbb{R}, \quad j_1 \simeq \alpha, \quad j_2(e^{i\tau \cdot}) \simeq T_{\tau} \quad \tau \in \mathbb{R}$$

Now we shall investigate the $C^*$-braided tensor product $A \otimes_{\mathcal{C}} C(\mathbb{T}^1)$.

For any $k \in C_{\mathrm{o}}(\mathcal{T}\mathbb{Z})$ we define the Fourier transform by:

$$(\mathcal{F}k)(p) := \frac{1}{\sqrt{2\pi}} \sum_n k(n\mathcal{T}) e^{-in\mathcal{T}p}$$

Notice that $\mathcal{F}k$ is a periodic function on $\mathbb{R}$ with the period $\frac{2\pi}{\mathcal{T}}: \mathcal{F}k \in C(\mathbb{T}^1)$.

In this case Formula (15) takes form

$$(\varphi v)(t) = \alpha_t(a) \sum_n k(n\mathcal{T}) v(t - n\mathcal{T}) =$$
$$= \alpha_t(a)\{\sum_n k(n\mathcal{T}) e^{-in\mathcal{T}\hat{p}} v\}(t) =$$
$$= \alpha_t(a)\{(\mathcal{F}k)(\hat{p})v\}(t)$$

Therefore

$$\varphi = [(\mathrm{id}_A \otimes i_1)\alpha(a)\{I_A \otimes i_2(\mathcal{F}k)\}$$

where $\mathcal{F}k$ runs over a dense subset of $C(\mathbb{T}^1)$ ([6] Theorem 1.2.4.).

Hence $(\mathrm{id}_A \otimes i_1)\alpha^A(A)(I_A \otimes i_2(C(\mathbb{T}^1)))$ is dense in $A \otimes_{\alpha} \mathcal{T}\mathbb{Z}$ (for simplicity of notation, we use te same letter $\alpha$ for action of the group $\mathbb{R}$ and $\mathcal{T}\mathbb{Z}$).

On the other hand ([8], Formula 2.5)

$$A \otimes_{\mathcal{C}} C(\mathbb{T}^1) \subset M\left(A \otimes C(\mathbb{T}^1) \otimes K\right) = C(\mathbb{T}^1, \mathrm{M}(A \otimes K))$$

where $C(\mathbb{T}^1, \mathrm{M}(A\otimes K))$ is a set of all strictly continuous $\mathrm{M}(A\otimes K)$ - valued functions on $\mathbb{T}^1$. For any $\xi \in A \otimes_{\mathcal{C}} C(\mathbb{T}^1)$ and $t \in \mathbb{R}$, let $\xi[t] \in \mathrm{M}(A\otimes K)$ denote the value of the function $\xi \in C(\mathbb{T}^1, \mathrm{M}(A \otimes K))$ at the point $t \in \mathbb{R}$:

$$\xi[t] = (\mathrm{id}_A \otimes \chi_t \otimes \mathrm{id}_K)\xi$$

where $\chi_t \in \mathrm{Mor}(C_\infty(\mathbb{R}), \mathbb{C})$ is the evaluation functional. With this notation definition (3) can be rewritten as follows

$$(j_1(a))[t] = (\mathrm{id} \otimes i_1)\alpha^A(a)$$

and

$$(j_2(b))[t] = (I \otimes (\mathrm{id} \otimes i_2)\alpha^B(b))[t] = I \otimes b(t+\hat{p}) = I \otimes e^{-it\hat{q}}b(\hat{p})e^{it\hat{q}} =$$
$$= (I \otimes e^{-it\hat{q}})(I \otimes i_2(b))(I \otimes e^{it\hat{q}})$$

Hence, observing that $j_1(a)$ and $e^{it\hat{q}}$ commute, we obtain

$$(j_1(a)j_2(b))[t] = (I \otimes e^{-it\hat{q}})(\mathrm{id} \otimes i_1)\alpha^A(a)(I \otimes i_2(b))(I \otimes e^{it\hat{q}})$$

and therefore $(\mathrm{id} \otimes i_1)\alpha^A(A)(I \otimes i_2(C(\mathbb{T}^1)))$ is unitarily equivalent to $j_1(A)j_2(C(\mathbb{T}^1))$, which is (see Definition 1.1) dense in $A \otimes_{\mathcal{C}} C(\mathbb{T}^1)$. Thus we see that $C^*$-algebras $A \otimes_{\mathcal{C}} C(\mathbb{T}^1)$ and $A \otimes_\alpha \mathcal{T}\mathbb{Z}$ are isomorphic

$$A \otimes_{\mathcal{C}} C(\mathbb{T}^1) \simeq A \otimes_\alpha \mathcal{T}\mathbb{Z}, \quad j_1 \simeq \alpha, \quad j_2(e^{i\tau \cdot}) \simeq T_\tau \quad \tau \in \mathcal{T}\mathbb{Z}$$

Now we shall investigate the $C^*$-braided tensor product $A \otimes_{\mathcal{C}} C(\mathbb{T}^2)$. In this case $G = \mathcal{T}_1\mathbb{Z} + \mathcal{T}_2\mathbb{Z} \subset \mathbb{R}$, where $\mathcal{T}_1, \mathcal{T}_2 \in \mathbb{R}_+$ and $\mathcal{T}_1/\mathcal{T}_2$ is not rational (the case when $\mathcal{T}_1/\mathcal{T}_2$ is rational is not interesting since it can be reduced to the previous example). Topology on $G$ is discret. For any $k \in C_o(\mathcal{T}_1\mathbb{Z} + \mathcal{T}_2\mathbb{Z})$ we define the Fourier transform by:

$$(\mathcal{F}k)(p) := \frac{1}{\sqrt{2\pi}} \sum_{n,m} k(n\mathcal{T}_1 + m\mathcal{T}_2)e^{-i(n\mathcal{T}_1+m\mathcal{T}_2)p}$$

Notice that $\mathcal{F}k$ is a finite linear combination of functions $e^{in\mathcal{T}_1 p + im\mathcal{T}_2 p}$, where $p \in \mathbb{R}$.

In this case Formula (15) takes form

$$(\varphi v)(t) = \alpha_t(a) \sum_{n,m} k(n\mathcal{T}_1 + m\mathcal{T}_2)v(t - n\mathcal{T}_1 - m\mathcal{T}_2) =$$
$$= \alpha_t(a)\{\sum_{n,m} k(n\mathcal{T}_1 + m\mathcal{T}_2)e^{-in\mathcal{T}_1\hat{p} - im\mathcal{T}_2\hat{p}}v\}(t) =$$
$$= \alpha_t(a)\{(\mathcal{F}k)(\hat{p})v\}(t)$$

Therefore

$$\varphi = [(\mathrm{id}_A \otimes i_1)\alpha(a)\{I_A \otimes i_2(\mathcal{F}k)\}$$

where $\mathcal{F}k$ runs over a set of finite linear combination of functions $e^{in\mathcal{T}_1 p + im\mathcal{T}_2 p}$, where $p \in \mathbb{R}$ ([6] Theorem 1.2.4.).

Observe that through a map $\rho$:

$$\rho:\ \mathbb{R} \ni t \mapsto (t + \frac{2\pi}{\mathcal{T}_1}\mathbb{Z}, t + \frac{2\pi}{\mathcal{T}_2}\mathbb{Z}) \in \mathbb{T}^2$$

the norm closure of a set: $\{\ \mathcal{F}k:\ k \in C_o(\mathcal{T}_1\mathbb{Z} + \mathcal{T}_2\mathbb{Z})\ \}$ is isomorphic to the set of all continuous functions on $\mathbb{T}^2$.

Hence $(\mathrm{id}_A \otimes i_1)\alpha^A(A)(I_A \otimes i_2(C(\mathbb{T}^2)))$ is dense in $A \otimes_\alpha (\mathcal{T}_1\mathbb{Z} + \mathcal{T}_2\mathbb{Z})$ (for simplicity of notation, we use te same letter $\alpha$ for action of the group $\mathbb{R}$ and $\mathcal{T}_1\mathbb{Z} + \mathcal{T}_2\mathbb{Z}$). On the other hand ([8], Formula 2.5)

$$A \otimes_\mathcal{C} C(\mathbb{T}^2) \subset M\left(A \otimes C(\mathbb{T}^2) \otimes K\right) = C(\mathbb{T}^2, \mathrm{M}(A \otimes K))$$

where $C(\mathbb{T}^2, \mathrm{M}(A \otimes K))$ is a set of all strictly continuous $\mathrm{M}(A \otimes K)$ - valued functions on $\mathbb{T}^2$, which we can treat also (using $\rho$) as functions on $\mathbb{R}$. For any $\xi \in A \otimes_\mathcal{C} C(\mathbb{T}^2)$ and $t \in \mathbb{R}$, let $\xi[t] \in \mathrm{M}(A \otimes K)$ denote the value of the function $\xi \in C(\mathbb{T}^2, \mathrm{M}(A \otimes K))$ at the point $t \in \mathbb{R}$:

$$\xi[t] = (\mathrm{id}_A \otimes \chi_{\rho(t)} \otimes \mathrm{id}_K)\xi$$

where $\chi_{\rho(t)} \in \mathrm{Mor}(C(\mathbb{T}^2), \mathbb{C})$ is the evaluation functional. With this notation definition (3) can be rewritten as follows

$$(j_1(a))[t] = (\mathrm{id} \otimes i_1)\alpha^A(a)$$

and

$$(j_2(b))[t] = (I \otimes (\mathrm{id} \otimes i_2)\alpha^B(b))[t] = I \otimes b(t + \hat{p}) = I \otimes e^{-it\hat{q}}b(\hat{p})e^{it\hat{q}} =$$
$$= (I \otimes e^{-it\hat{q}})(I \otimes i_2(b))(I \otimes e^{it\hat{q}})$$

Hence, observing that $j_1(a)$ and $e^{it\hat{q}}$ commute, we obtain

$$(j_1(a)j_2(b))[t] = (I \otimes e^{-it\hat{q}})(\mathrm{id} \otimes i_1)\alpha^A(a)(I \otimes i_2(b))(I \otimes e^{it\hat{q}})$$

and therefore $(\mathrm{id} \otimes i_1)\alpha^A(A)(I \otimes i_2(C(\mathbb{T}^2)))$ is unitarily equivalent to $j_1(A)j_2(C(\mathbb{T}^2))$. which is (see Definition 1.1) dense in $A \otimes_\mathcal{C} C(\mathbb{T}^2)$.

Thus we see that $C^*$-algebras $A \otimes_\mathcal{C} C(\mathbb{T}^2)$ and $A \otimes_\alpha (\mathcal{T}_1\mathbb{Z} + \mathcal{T}_2\mathbb{Z})$ are isomorphic

$$A \otimes_\mathcal{C} C(\mathbb{T}^2) \simeq A \otimes_\alpha (\mathcal{T}_1\mathbb{Z} + \mathcal{T}_2\mathbb{Z}), \quad j_1 \simeq \alpha, \quad j_2(e^{i\tau\cdot}) \simeq T_\tau \quad \tau \in \mathcal{T}_1\mathbb{Z} + \mathcal{T}_2\mathbb{Z}$$

### Acknowledgements

This is an extended version of the author's M. Sc. thesis, written under the supervision of Prof. S. L. Woronowicz at the University of Warsaw. The author is greatly indebted to Prof. S. L. Woronowicz for his assistance during the preparation of the paper.

# References

[1] O. Bratelli, D. Robinson: *Operator Algebras and Quantum Statistical Mechanics, Part I*, Springer-Verlag, New York, Heidelberg, Berlin 1979.

[2] M.B.Landstad: *Transactions of the American Mathematical Society*, **248**, 223 (1979).

[3] S.MacLane: *Categories for Working Mathematician*, Springer-Verlag, New York 1971.

[4] S.Majid: *Foundations of Quantum Group Theory*, Cambridge University Press, Cambridge 1995.

[5] G.K. Pedersen: $C^*$ - *Algebras and their Automorphism Groups*, Academic Press London, New York, San Francisco 1962.

[6] W. Rudin: *Fourier Analysis on Groups*, Interscience Publishers, New York, London 1962.

[7] S.L. Woronowicz: *An example of a braided locally compact group* in *The Proc. of XXX–th Karpacz Winter School of Theor. Phys.*, PWN, Warszawa 1995.

[8] S.L. Woronowicz: *Reviews on Mathematical Physics*, **7**, 481 (1995).